\documentclass{article}
\usepackage[normalem]{ulem}
\usepackage{amsthm}
\usepackage{mathtools}
\usepackage{thmtools}

\declaretheoremstyle[headfont=\normalfont]{normalhead}
\usepackage{color}
\usepackage{graphicx}
\usepackage{morefloats}
\usepackage{hyperref}
\usepackage{mathrsfs}
\hypersetup{
    colorlinks=false,
    linkcolor=cyan,
    filecolor=cyan,      
    urlcolor=cyan,
}
\usepackage{biblatex}
\usepackage{ amssymb }
\usepackage{textcomp}
\usepackage{mathabx }
\usepackage{url}
\usepackage{float}
\usepackage{biblatex}
\newtheoremstyle{mydef}
{\topsep}{\topsep}%
{}{}%
{\itshape}{}
{\newline}
{%
  \rule{\textwidth}{0.0pt}\\*%
  \thmname{#1}~\thmnumber{#2}\thmnote{\-\ #3}.\\*[-1.5ex]%
  \rule{\textwidth}{0.0pt}}%

\begin{document}
\newcommand{\rb}[1]{\raisebox{-.2ex}[0pt]{#1}}
\newtheorem{conjecture}{Conjecture}
\newtheorem{theorem}{Theorem}
\newtheorem*{theorem-non}{Theorem}
\newtheorem{remark}{Remark}
\newtheorem{proposal}{Proposal}
\newtheorem{proposition}{Proposition}
\newtheorem*{proposition-non}{Proposition}
\newtheorem{lemma}{Lemma}
\newtheorem{corollary}{Corollary}
\newtheorem{observation}{Observation}
\newtheorem{definition}{Definition}
\newtheorem*{question-non}{Question}
\newtheorem*{corollary-non}{Corollary}
\author{Barry Brent}

\date{23 December 2022}

\title{On the constant terms of certain meromorphic modular
forms for Hecke groups}
\maketitle
\begin{abstract}\noindent
We study polynomials interpolating
the (rational) constant terms of 
certain meromorphic modular forms
for Hecke groups. We make observations
about the divisibility properties
of the constant terms and connect them
to several sequences, for example, 
to O.E.I.S. sequence A005148 
\cite{OEISNewmanShanks},
which was studied by Newman, Shanks
and Zagier
\cite{newman2004sequence}, 
\cite{newman2004sequenceAppendix}
in an article on its use in series 
approximations to $\pi$.
\hskip -.2in
\end{abstract}
\section{Introduction}
The study of the 
constant terms of meromorphic modular forms
bears upon the analysis  of ordinary 
quadratic forms.
C. L. Siegel studied the constant terms in the
Fourier expansions of a particular family 
of meromorphic modular forms $T_h$ for
$SL(2,\mathbb{Z})$ (``level one modular forms'') 
in 1969
\cite{siegel1969berechnung, siegel1980evaluation}.
In the relevant 
part of his article, Siegel
demonstrated that these constant terms never vanish.
He used this to 
establish a bound on 
the exponent of the first non-vanishing Fourier 
coefficient for
a level one entire modular form $f$ of weight 
$h$ such that the
constant term of $f$ is itself non-vanishing. 
Theta 
functions fit this description, so Siegel 
was able to give an upper 
bound on the least positive integer 
represented by a positive-definite
even unimodular quadratic form in $2h$ variables.
\newline \newline \noindent
While looking at the level two situation,
the present writer found 
numerical evidence for
divisibility properties of the constant terms 
for several kinds of modular form, including 
the $T_h$ \cite{brent1998quadratic}; 
if these properties hold, the constant terms 
cannot vanish. 
To conform to our notation in the sequel,
let $A_{\overline{\mathcal{K}},k,3}(0)$ 
be the constant
term of $j^k$ where $j$ is the usual  
Klein invariant. 
(Thus $A_{\overline{\mathcal{K}},1,3}(0) = 744$;
we defer our explanation of the symbol 
$\overline{\mathcal{K}}$.)
\footnote{For example, see Serre \cite{serre1970course}, 
section 3.3, equation (22), or the Wikipedia page
\cite{jWiki}.}
Furthermore let $d_b(n)$ be the sum of the digits in 
the base $b$
expansion of $n$. Then (apparently) 
$$\text{ord}_2(A_{\overline{\mathcal{K}},k,3}(0)) 
= 3d_2(k)$$
and $$\text{ord}_3(A_{\overline{\mathcal{K}},k,3}(0))) = 
d_3(k).$$ 
In this 
article we will argue, but only
empirically, that the 
$A_{\overline{\mathcal{K}},k,3}(0)$ 
inherit the stated properties from the 
OEIS sequence A005148 \cite{OEISNewmanShanks},
which was originally studied by Newman, Shanks
and Zagier
\cite{newman2004sequence, 
newman2004sequenceAppendix}
in an article on its use in series 
approximations to $\pi$. The constant terms
in the Fourier expansions of other modular
forms appear to inherit such divisibility 
properties from sequences described below
that (so far) are not included in 
Sloane's encyclopedia.
\newline \newline \noindent
We tried 
to find patterns in the
$p$-orders of constant terms of 
$j$ and other modular forms for
$SL(2,\mathbb{Z})$ for $p$ larger than three.
When our search failed, we 
began to search among the  Hecke groups,
because $SL(2,\mathbb{Z})$ is the first one, 
isomorphic to the product of
cyclic groups $C_2*C_3$,
while in general they 
have the form $C_2*C_m$ for 
for $m = 3, 4, ....$
We will state some conjectures
about the constant terms, for example, of
meromorphic forms for
Hecke groups isomorphic to $C_2*C_{p^k}, p$
prime.
\section{Background}
For $m = 3, 4, ...$,
let $\lambda_m = 2 \cos \pi/m$ and
let $J_m$ be a certain meromorphic
modular form for the Hecke group $G(\lambda_m)$,
built from triangle functions, 
with Fourier expansion
$$
J_m(\tau) = \sum_{n=-1}^\infty a_n(m)q_m^n,
$$
where
$q_m(\tau) = \exp 2 \pi i \tau/\lambda_m$.
The groups $G(\lambda_m)$ and $
SL(2,\mathbb{Z})$ coincide.
(For further details, the reader
is referred to the books by 
 Carath{\'e}odory
\cite{caratheodory1, caratheodory2}
and by Berndt and Knopp
\cite{berndt2008hecke},
the articles of Lehner and Raleigh
\cite{lehner1954note, raleigh1962fourier},
to the dissertation of Leo 
\cite{leo2008fourier}, and to a
summary, including pertinent references
to that material, in the
2021 article \cite{interpolating}.)
\newline \newline \noindent
Raleigh gave polynomials $P_n(x)$ such that 
$a_{-1}(m)^n q_m^{2n+2} a_n(m) = P_n(m)$
for  $n = -1, 0, 1, 2$ and $3$.
He conjectured that similar 
relations hold for all positive integers
$n$ \cite{raleigh1962fourier}.
\footnote{For more on
expansions over polynomial fields, see, for 
example, the book of Boas and Buck 
\cite{boas2013polynomial} and the articles 
by Buckholtz and Byrd (\cite{buckholtz1973series},
\cite{byrd}.)}
Akiyama proved Raleigh's conjectures
in 1992 \cite{akiyama1992note}.
\newline \newline \noindent
Using the weight-raising properties of 
differentiation and the $J_m$, Erich
Hecke constructed
certain families 
$\mathcal{H}$
comprising modular forms 
of positive weight
for each $G(\lambda_m)$
sharing certain properties 
\cite{hecke1936, berndt2008hecke}. 
(The weight of $g$ is not necessarily
constant within such a family.)
It seems apparent that Akiyama's 
result can be extended:
there should exist polynomials 
$Q_{\mathcal{H},n}(x)$ 
interpolating the coefficient of $X_m^n$
in the Fourier expansions of the
members of Hecke families
$\mathcal{H}$. \footnote{An
incomplete draft
of a longer
article justifying this idea is
located in the folder `current draft' 
\cite{githubNewmanShanks}.
but we have already acted on it in our 2021 
Integers paper \cite{interpolating}.}
\newline \newline \noindent
In section 4 of our 2021 article,
we made use of a certain uniformizing
variable $X_m(\tau)$ 
for $\tau$ in the upper half 
plane \cite{interpolating}.
By Akiyama's theorem, we have a
series of the form
$\mathcal{J}_m(x) := 
\sum_{n=-1}^{\infty} \tilde{P}_n(x)X_m^n$
for polynomials $\tilde{P}_n(x)$ 
in $\mathbb{Q}[x]$
with the property that 
$J_m = \mathcal{J}_m(m)$.
We will make use of the change of
variables  $X_m \mapsto 2^6 m^3 X_m$ for a 
$G(\lambda_m)$-modular form (originally employed,
as far as we know, by Leo (\cite{leo2008fourier},
page 31). It has the effect when $m=3$
of recovering the Fourier series
of a variety of standard 
modular forms. We set this up as a 
\begin{definition}
For $\tau$ in the half plane
$\{z \in \mathbb{C}$ such that $\Im(z) >0\}$
\footnote{This is the usual domain of a
classical modular form or modular function.}
and $k_a \neq 0$, let 
$$
f(\tau) = \sum_{n=a}^{\infty} k_n X_m(\tau)^n
$$
and
$$
g(\tau) = \sum_{n=a}^{\infty} k_n 2^{6n} m^{3n} 
X_m(\tau)^n.
$$
If we rewrite the last expansion as 
$g(\tau)  = \sum_{n=a}^{\infty} \tilde{k}_n X_m(\tau)^n$,
then we set
$$\overline{f}(\tau) := g(\tau)/\tilde{k}_a.$$
Also, for $m  = 3, 4, ...$,
we set $j_m(\tau) := \overline{J_m}(\tau)$.
\end{definition} \noindent
The Fourier expansion of $j_3$
is 
\footnote{
See equation (23) of 
Serre's book 
\cite{serre1970course}, 
section 3, and the \it SageMath \rm notebook 
``jpower constant term NewmanShanks 26oct22.ipynb''
in 
\cite{githubNewmanShanks}.
}
$$j_3(\tau) = 
1/X_3(\tau) + 744 + 196884 X_3(\tau) + 
21493760 X_3(\tau)^2 + ...,
$$
which matches the  standard expansion $j(\tau) = $
$$
1/\exp(2\pi i \tau) + 744 + 
196884 \exp(2\pi i \cdot  \tau) + 
21493760 \exp(2\pi i\cdot  2 \cdot \tau ) + ....
$$
 \section{Fourier expansions}
 \noindent
 We make the following 
 \begin{definition}
 Let $\mathcal{F} = \{f_3, ..., f_m, ...\}$
 where $f_m$ is modular for $G(\lambda_m)$.
 Then we write the Fourier expansion
 of $f_m^k$ in powers of $X_m$ as
$$
f_m(\tau)^k = \sum_n A_{\mathcal{F},k,m}(n)X_m^n.
$$
 \end{definition} \noindent
\begin{proposition}
Let 
$\mathcal{K} =
\{J_3, , J_4,...\}$
and $\overline{\mathcal{K}} =
\{j_3, j_4,....\}$
Then there exist polynomials
$Q_{\mathcal{K},k,n}(x)$ and
$Q_{\overline{\mathcal{K}},k,n}(x)$
in $\mathbb{Q}[x]$
such that 
$$J_m(\tau)^k = 
\sum_{n = -k}^{\infty} 
Q_{\mathcal{K},k,n}(m) X_m(\tau)^n$$ and
$$j_m(\tau)^k = 
\sum_{n = -k}^{\infty} 
Q_{\overline{\mathcal{K}},k,n}(m) 
X_m(\tau)^n.$$
In other words,
$A_{\mathcal{K},k,m}(n) = Q_{\mathcal{K},k,n}(m)$
and 
$A_{\overline{\mathcal{K}},k,m}(n) = 
Q_{\overline{\mathcal{K}},k,n}(m)$
for $k = 1, 2, ..., m = 3, 4, ...$,
and $n = -k, 1-k, ....$
\end{proposition} 
\noindent
For $k$ equal to one, 
the first claim is
just Akiyama's theorem and
the claim for $k$ not equal to
one is then obvious.
The second statement 
follows immediately.
\newline \newline \noindent
When, given  a sequence of functions
$f_m$ modular for  $G(\lambda_m)$ 
in a family $\mathcal{F}$,
we wanted to find polynomials
$Q_{\mathcal{F},n}(x)$
such that each $f_m$
with Fourier expansion
$$
f_m(\tau) = \sum_n  a_{m,n} X_m^n(\tau),
$$
satisfied $Q_{\mathcal{F},n}(m) = a_{m,n}$,
we evaluated
finite sequences
$\{a_{m,n}\}_{m = 3, 4, ...,M}$
 (with $n$ held constant)
and generated the candidates
for $Q_{\mathcal{F},n}(x)$
by Lagrange interpolation.
The bound $M$ was 
chosen large enough that
the degrees of the $g_n(x)$
that the procedure produced were
linear in $n$.
Over the
course of experiments described in
our earlier article \cite{interpolating}, 
this linearity 
was associated with systematic behavior.
For example, if a polynomial $g_n(x)$
was factored as 
$g_n(x) = 
r_n\cdot p_1(x) \cdot p_2(x)...\cdot p_a(x)$
where each of the $p_i$ was monic, $r_n$
was rational, and the degree of $g_n(x)$
was linear in $n$, then often the sequence
$\{r_3, r_4, ...\}$ was readily identifiable
(sometimes after resorting to Sloane's
encyclopedia.) We take such regularities
as evidence that the polynomial $g_n(m) =
a_{m,n}$ for all $m$.
\section{Divisibility properties of constant terms
for weight zero meromorphic modular functions}
Let $\mathcal{K}, \overline{\mathcal{K}},
A_{\mathcal{K},k,m}(n)$,
and 
$A_{\overline{\mathcal{K}},k,m}(n)$ 
be as in proposition 1.
 When checking the conjectures in this 
 section\footnote{
See the \it SageMath \rm notebooks
in our repository
\cite{githubNewmanShanks},
in the folder ``renumbered conjectures''.
}, we did not
use  tables of the
$A_{\mathcal{K},k,m}(0)$ 
and $A_{\overline{\mathcal{K}},k,m}(0)$
directly.
Instead (for example), we used 
Lagrange interpolation
to identify polynomials $h_k(x)$
and $\overline{h}_k(x)$
such that $A_{\mathcal{K},k,m}(0) = h_k(m)$
and
$A_{\overline{\mathcal{K}},k,m}(0) = 
\overline{h}_k(m)$
by letting $m$ run through a small set
of values sufficient to produce
the linearity behavior we mentioned
in the previous section; thus
we have assumed (in this example) that
$h_k(x) \equiv Q_{\mathcal{K},k,0}(x)$
and
$\overline{h}_k(x) \equiv 
h_k(x) \equiv Q_{\mathcal{K},k,0}(x)$
and
$\overline{h}_k(x) \equiv 
Q_{\overline{\mathcal{K}},k,0}(x)$
identically. We made tables
of $p$ orders of the 
$h_k(m)$ and the $\overline{h}_k(m)$. 
In this way we checked larger sets of $m$
values
than would have been practicable if we 
had checked
the constant terms themselves.
Similar remarks will apply to our methods of
studying the
constant terms of negative weight meromorphic
modular forms: reciprocal powers of the cusp 
forms defined above.
\newline \newline \noindent
There are a variety of logical relations
among the conjectures below. We of course
do not know which of them (if any) are correct; 
so, for now, we state them separately.
\newline \newline \noindent
Unlike the later conjectures, conjecture 1
is not a way of summarizing patterns we saw
in our data. Rather it codifies our fundamental
assumption that the linearity behavior
we described is a reliable signal.
\begin{conjecture}
\begin{enumerate}
\item 
$h_k(x) \equiv 
Q_{\mathcal{K},k,0}(x)$ identically;
consequently, 
$h_k(m) = A_{\mathcal{K},k,m}(0)$
identically.
\item
$\overline{h}_k(x) \equiv 
Q_{\overline{\mathcal{K}},k,0}(x)$
identically; consequently,
$\overline{h}_k(m) = 
A_{\overline{\mathcal{K}},k,m}(0)$
identically.
\end{enumerate}
\end{conjecture} \noindent
\subsection[]{The constant terms 
$A_{\overline{\mathcal{K}},k,m}(0)$.
} 
\subsubsection[]{$m = 3$.}
In this subsection, our goal is to illustrate
connections between the divisibility
patterns (described in the
introduction)
for the constant terms of the
$j(\tau)^k$ Fourier expansions 
on one side,
and the $h_k(x)$ on the other.
Let $\overline{h}_k(x)$ factor as 
$\overline{h}_k(x) = 
\nu_k\cdot p_{k,1}(x) \times p_{k,2}(x) \times ...
\times p_{k,\alpha}(x) =$
(say)
$\nu_k\cdot \tilde{p}_k(x)$
where each of the $p_{k,n} (n = 1, 2,...,
\alpha$) is monic and $\nu_k$
is rational. 
Let us represent O.E.I.S.
sequence A005148 \cite{OEISNewmanShanks}  
$\{0, 1, 47, 2488, 138799, ...\}$
as $\{a_0, a_1, ...\}$.
\begin{conjecture}
\begin{enumerate}
\item
$\nu_k = 24a_k$.
\item $\tilde{p}_k(3)$
is always odd.
\item
$\text{ord}_2(a_k) = 3d_2(k)-3$.
\item
$\text{ord}_3(\tilde{p}_k(3)) = 
d_3(k) - 1$.
\end{enumerate}
\end{conjecture} \noindent
We have arrived at the observations from
our 1998 article as described in the introduction.
\begin{corollary}
$\text{ord}_2(A_{\overline{\mathcal{K}},k,3}(0)) = 
3d_2(k)$
and
$\text{ord}_3(A_{\overline{\mathcal{K}},k,3}(0)) = d_3(k)$.
\end{corollary} \noindent
\it{Proof} \rm \thinspace
First claim:
$\text{ord}_2(A_{\overline{\mathcal{K}},k,3}(0)) =
\text{ord}_2(\overline{h}_k(3)) =
\text{ord}_2(\nu_k \cdot \tilde{p}_k(3)) =
\text{ord}_2(24 a_k) \cdot \tilde{p}_k(3)) =
\text{ord}_2(24) + \text{ord}_2(a_k)+ 
\text{ord}_2(\tilde{p}_k(3)) =
3 + 3d_2(k)-3 + 0 = 3d_2(k).$
 Second claim:
In their 1984 article 
\cite{newman2004sequence},
Newman, Shanks and Zagier demonstrated
that
$\text{ord}_3(a_k) = 0$
for all $k$.
Therefore (under the conjectures above)
$\text{ord}_3(A_{\overline{\mathcal{K}},k,3}(0))
=\text{ord}_3(\overline{h}_k(3)) =
\text{ord}_3(\nu_k) +
\text{ord}_3(\tilde{p}_k(3)) =
1 + \text{ord}_3(a_k) + d_3(k) - 1 =
d_3(k)$.
\subsubsection[]{$m$ a prime power.}
By imposing 
restrictions on $k$ and $m$,
we found several narrow
conjectures  about constant term $p$ orders 
for various primes $p$.
\footnote{Again, see
the \it SageMath \rm
notebooks
in the folder ``renumbered conjectures''
in our repository
\cite{githubNewmanShanks}.
We identified the sequences 
involved after reading several pages in the O.E.I.S. 
\cite{OEISturpel},\cite{OEISgerasimov},
\cite{OEISsloane},\cite{OEISzum}.}
\begin{conjecture}
If $p$ is prime
and $a$ is an integer that is 
larger than $2$,
then
$$\text{ord}_p(A_{\overline{\mathcal{K}},k,p^a}(0)) 
=(a-3)k+
\text{ord}_p(A_{\overline{\mathcal{K}},k,p^3}(0))).
$$
\end{conjecture}
\begin{conjecture}
Let $a \geq 2$. Then
$\text{ord}_2(A_{\overline{\mathcal{K}},
2,2^a}(0)) = 2a+7$.
\end{conjecture}
\begin{conjecture}
Let $p$ be a prime number larger than $2$
and let $a$ be a positive integer. Then
$\text{ord}_p(A_{\overline{\mathcal{K}},p,p^a}(0))
= ap-2$.
\end{conjecture} \noindent
\subsubsection[]{Other $m$.}
\begin{conjecture}
If $d_2(k) = 1$, 
$a = \text{ord}_2(m), 
a \geq 2$, and
$o = \text{ord}_2 (A_{\overline{\mathcal{K}},k,m}(0))$,
then $o=k(a+2)+3$.
\end{conjecture}
\begin{conjecture} 
Let $d_2(k) = 1$,
$m \equiv 2 \thinspace 
(\text{mod} \thinspace 4)$,
and $a=\text{ord}_2(m) ( = 1$,
of course.)
Then 
$\text{ord}_2 (A_{\overline{\mathcal{K}},k,m}(0)) =
k(a+6) + 1 = 7k + 1$.
\end{conjecture} \noindent
Now let $C_n, n = 0, 1, 2, ...$ be 
the $n^{th}$ Catalan number. 
One of several explicit formulas 
for $C_n$ is $$C_n = 
\frac{(2n)!}{(n+1)!n!}.$$ \noindent
For $n$ positive
let  $C_{1,n}$ denote the $n^{th}$ 
Catalan number $c$ such that $c \neq C_0$ and $\text{ord}_2(c) = 1$.\footnote{We 
 encountered 
 it on Bottomley's O.E.I.S. page 
\cite{OEISbottomley}.}
\begin{conjecture}
Let 
$k$ be the $n^{th}$ positive integer
such that
$d_2(k) = 2$; also,
$m = 4j, (j = 1, 2,...)$,
and
$a = \text{ord}_2(m)$.
Furthermore, let
$o = \text{ord}_2 (A_{\overline{\mathcal{K}},k,m}(0)$ and 
$ t = ((a+6)k +2 - o)/4$.
Then 
$t = C_{1,n}$.
\end{conjecture}
\begin{conjecture}
Let  $d_2(k) = 2$,
$m = 4j + 2, j = 1, 2, ...$, and 
$a = \text{ord}_2(m)$ (again: $ = 1$.)
Then
$\text{ord}_2 
(A_{\overline{\mathcal{K}},k,m}(0)) = (a+6)k+2
= 7k + 2$.
\end{conjecture}
\begin{conjecture} If
$m \equiv 0 \thinspace 
(\text{mod} \thinspace 3)$,
then $\text{ord}_3 (A_{\overline{\mathcal{K}},k,m}(0))=
k \cdot \text{ord}_3(m) + d_3(k) - k.$
\end{conjecture}
\subsection[]{The constant terms $A_{\mathcal{K},k,m}(0)$.
} 
\begin{conjecture}
\begin{enumerate}
\item
 Let $p$ be a prime number greater than two.
Then\footnote{The supporting data files
for this conjecture are too large to store
conveniently on GitHub; instead, we stored
links to them there in a file
called `links'.\cite{githubNewmanShanks} (The links download the files from the author's Google 
drive.)}
$$\text{ord}_p(A_{\mathcal{K},p,p}(0)) = 
-2 - 2p.$$
\item The value of $\text{ord}_p(A_{\mathcal{K},p,p^n}(0))$
does not vary on the set $n = 2, 3, ...$.
\item  If $n$ is an integer greater than one
and $d_2(p)>2$, 
\footnote{This version of the
restriction on
$d_2(p)$ is supported by only a little evidence, since exceptions to it
are scarce in our data set.
See \cite{OEIScoveiro} and
the file `conjecture 11 no2 23dec22.ipynb'
in the folder `renumbered conjectures'
in our depository
\cite{githubNewmanShanks}.}
$$\text{ord}_p(A_{\mathcal{K},p,p^n}(0)) = 
-2 - 2p.$$
\end{enumerate}
\end{conjecture}
\printbibliography
$\tt{barrybrent@iphouse.com}$
\end{document}